\documentclass{article}
\usepackage{graphicx}
 \usepackage{mathptmx}
\usepackage{amsmath, amstext,amssymb,amsfonts}
\usepackage[english]{babel}
\usepackage{delarray}
\usepackage{mathptmx}
\usepackage{amsmath, amstext,amssymb,amsfonts}
\usepackage[english]{babel}
\usepackage{delarray}

\DeclareFontFamily{U}{mathx}{\hyphenchar\font45}
\DeclareFontShape{U}{mathx}{m}{n}{
      <5> <6> <7> <8> <9> <10>
      <10.95> <12> <14.4> <17.28> <20.74> <24.88>
      mathx10
      }{}
\DeclareSymbolFont{mathx}{U}{mathx}{m}{n}
\DeclareFontSubstitution{U}{mathx}{m}{n}
\DeclareMathAccent{\widecheck}{0}{mathx}{"71}
\DeclareMathAccent{\wideparen}{0}{mathx}{"75}

\numberwithin{equation}{section}
  \newcommand{\be}{\,B^1_{\sigma}\,}  \newcommand{\bep}{\,B^{p}_{\sigma}\,}
 \newcommand{\R}{\mathbb R}   \newcommand{\Co}{\mathbb C}
 
\newcommand{\ep}{\varepsilon}
 \newcommand{\dr}{{\text d}}
 
 \newcommand{\dbe}{\ {\cal{D}}(B^1_{\sigma})\,}        
 \newcommand{\expp}{   {\text{ exp}}\dbe}                 \newcommand{\exps}{{\text{ exps}}\dbe}
  \newcommand{\he}{H^{\infty}_{\sigma}(\R) }
 \newcommand{\exppp}{   {\text{ exp}}^{\ast}\dbe}             \newcommand{\expss}{{\text{ exps}}^{\ast}\dbe}

\newtheorem{teo}{Theorem}[section]    

\date{}
\title{On exposed functions in Bernstein spaces of functions of exponential type}
\author{\large{ Saulius  Norvidas}}
\date{\footnotesize Institute of Data Science and Digital Technologies, Vilnius University, \\ Akademijos str. 4, Vilnius LT-04812, Lithuania\\
 ({\rm{e-mail: norvidas{@}gmail.com}})}
\begin{document}

\maketitle
 {{ {\bf Abstract}}}  \ \ For  $\sigma>0$, the Bernstein space \ $B^1_{\sigma}$ consists of those  $L^1(R)$\
functions whose Fourier transforms are supported by $[-\sigma,\sigma]$.  Since $B^1_{\sigma}$ is separable and dual to
some Banach space,   the closed unit ball  $D(B^1_{\sigma})$ of $B^1_{\sigma}$\  has sufficiently large sets of both exposed and
strongly exposed points. Moreover,   $D(B^1_{\sigma})$  coincides with the closed convex  hull of its strongly exposed points. We investigate
some properties of exposed points, construct several examples and obtain as corollaries the relations between the
sets of exposed, strongly exposed, weak$^{\ast}$  exposed, and weak$^{\ast}$  strongly exposed points of  $D(B^1_{\sigma})$.



\section{ Introduction }
{\large{
In this paper we study the structure  of the set of exposed points in the unit ball of certain proper
subspace of \ $L^{1}(\R)$.\  For\ $0<\sigma<\infty$,\ the closed subspace we consider is the Bernstein space\
$\be$\ of complex valued\ $[-\sigma, \sigma]$--bandlimited functions in\ $L^1(\R)$,\ i.e.,
 \[
 \be=\{f\in L^1(\R): \ {\text{supp}}\,{\hat{f}}\subset [-\sigma, \sigma]\},
 \]
 where\ $ {\hat{f}}$\ is the Fourier transform of\ $f$.\ The class $\be$\  is the Banach spaces when endowed with the natural\ $  L^1(\R)$-norm.
 Denote by\  $\dbe$\  the closed unit ball of\ $\be$.\ Recall
 that\ $f\in\dbe$\ is  {\it extreme} function (point) of\ $\dbe$\  if for any \ $u, v\in\dbe$,\ $f=\frac12(u+v)$,\
 implies that\ $u=v=f$.\ Next,  a function\ $f\in\dbe$\  is said to be {\it exposed}   if there exists a linear functional\
 $\Phi_f$\ on\ $\be$\ such that $\Phi_f(f)=\|\Phi_f\|=1$\ and\ ${\text{Re}}\, \Phi_f(g)<1$\ for
 all\ $g\in\dbe$,\ $g\neq f$.\  The functional\ $\Phi_f$\ is  called an  {\it exposing functional} for\ $f$.   Finally,
 we  call\ $f\in\expp$\
 {\it strongly exposed} if for every sequence\ $(f_n)$\ in\ $\dbe$\ such that\ $\Phi_f(f_n)\to \Phi_f(f)=1$,\ it
 follows that\ $f_n\to f$\ in\ $\be$.\  We write\ $\expp$\ for the set  of exposed points in\ $\dbe$,\  and the set
 of strongly exposed  points in\ $\dbe$\ we denote by\ $\exps$.\

 \ \ \ \ \  The  space\ $\be$\  is dual to some Banach space (the details are contained in the next section). Thus, the
 following Phelps theorem guarantees  that\ $\exps$\   are sufficiently
 large subsets  of\ $\dbe$\ [8]: in a separable dual Banach space the closed unit ball coincides with the closed convex
 hull of its strongly exposed points. The extreme points of\ $\dbe$\ have a simple characterization in terms of
 zeros  of entire functions, whereas description of exposedness is less transparent (see Theorem 2.1). We also note
 that strongly exposed points have not been studied so far.

 \ \ \ \ \ Given\ $f\in\expp$,\ its exposing functional\ $\Phi_f$\ is determined uniquely (see  section 3). Therefore,
 characterization of\ $\expp$\ is  equivalent to  description of functionals\
 $\Psi\in(\be)^{\ast}$,\ $\|\Psi\|=1$,\  which attain its norm and have a
 unique {\it extremal},  i.e., there exists a unique\ $f\in\dbe$\ such that\ $\Psi(f)=\|\Psi\|=1$.
 We select weak$^{\ast}$\ continuous  functionals on\ $\be$\ as the most simple functionals, which attain its norm.
     Let\ $f\in\expp$,\ and suppose that exposing functional\ $\Phi_f$\  is weak$^{\ast}$\ continuous.
 Such a function\ $f$\ is said to  be {\it exposed}$^{\ast}$\ (or weak$^{\ast}$\ exposed). The set of exposed$^{\ast}$\
 functions  in\ $\dbe$\ will be denoted by\ $\exppp$.\  We call\ $f\in\exppp\cap\exps$\ {\it strongly exposed}$^{\ast}$\
 (or weak$^{\ast}$ strongly exposed), and we denote the set\ $ \exppp\cap\exps$\ by\ $\expss$.\ It is known\ [1]\ the
 following result attributed to E. Bishop: if\ $E$\ is a Banach space with separable dual\ $E^{\ast}$\ and\ $K$\
 is a  non-empty convex weak$^{\ast}$ compact subset of\ $E^{\ast}$,\ then\ $K$\ is the weak$^{\ast}$ closed convex
 hull of\  ${\text{exps}}^{\ast}(K)$.

\ \ \ \ \ In the present paper,  relations between\ $\expp$,\ $\exps$,\
$\exppp$,\ and\ $\expss$\ are studied. It is obvious that\ $\expss\subset \exppp\subset\expp$,\ and\
 $\expss\subset \exps\subset \expp$.\  We prove some properties of exposed functions, construct several examples, and
 obtain as corollaries the following  new relations
\[
\expss=\exppp,
\]
and
\[
\exppp\subsetneqq\exps\subsetneqq\expp.
\]
\vspace{4mm}
\section{  PRELIMINARIES }

\ \ \ \ \ To make the article self-contained, we recall   some definitions and facts about entire
functions and geometry of\ $\dbe$.\  An entire function\ $f:\mathbb{C}\to\mathbb{C}$\  is said to be of {\it exponential
type}  at most\ $a$ \ ($0\le a <\infty)$\  if for every\ $\ep >0$\ there  exists an\ $M_{\ep}>0$\ such that
\[
|f(z)|\leqslant  M_{\ep}e^{(a+\ep)|z|},\quad z\in\Co.
\]
The greatest lower bound of those\ $a$\ coincides with the {\it type\ $\sigma_f$\ of entire function}\ $f$\ (concerning
the first order growth). For\ $1\le p\le\infty$\ and\ $0<\sigma<\infty$,\ the Bernstein space\ $B_{\sigma}^p$\ is defined
as the set of functions\ $f\in L^p(\R)$\ such that their Fourier transforms (in the distributional sense) vanish outside\
$[-\sigma, \sigma]$.\ These\ $B_{\sigma}^p$\  are Banach spaces with respect to the\ \ $ L^p(\R)$--norm, and are ordered
by imbending\  $\be\subset B_{\sigma}^p \subset B^r_{\sigma}\subset B^{\infty}_{\sigma}$, $1\leqslant p\leqslant r
\leqslant {\infty}$,\ [3].
By the Paley-Wiener theorem and its inversion, the space\ $B_{\sigma}^p$\ can be described as the class
of\ $f\in L^P(\R)$\ such that\ $f$\ has  an extension from\ $\R$\ onto the complex plane\ $\Co$\  to an entire function of
exponential type at most\ $\sigma$.\ Therefore, we can identify\ $f\in B_{\sigma}^p\subset L^{p}(\R)$\ with its
entire extension on\ $\Co$.\

\ \ \ \ \ Every\ $f\in B_{\sigma}^p$,\ $1\le p<\infty$,\  satisfies the Plancherel--Polya inequality
\begin{equation}
\int\limits_{-\infty}^{\infty}|f(x+iy)|^p\,d x\le e^{p \sigma|y|}\int\limits_{-\infty}^{\infty}|f(x)|^{p}\,d x=
e^{p \sigma|y|}\|f\|_{L^{p}(\R)}^{p},\quad y\in\R,
\end{equation}
and the condition\ [3]
\begin{equation}
\lim_{x\to\pm\infty} f(x)= 0.
\end{equation}
\ \ \ \ From\ (2.2)\ it follows that any nonzero\ $f\in \bep$,\ $1\le p\le\infty$,\ has in\ $\Co$\ infinitely many
zeros (roots). Denote by\ $N(f)$\ the set of all zeros of\ $f\in\bep$\  with multiplicities  counted. If\
$\lambda\in N(f)\cap\R$,\  we call such\ $\lambda$\
  a {\it real zero} of\ $f$,\ and we use the term  a {\it complex zero} only for\   $\lambda\in N(f)$\ belonging
  to\ $ \Co\setminus\R$.\   If\  $f(\lambda)=f({\overline \lambda})=0$,\  we call
such\ $\lambda\in N(f)$\  a {\it  conjugate  zero} of\ $f$.\  For\ $f\in\bep$,\  we define its  {\it conjugate
function} as\  $f^{\ast}(z):={\overline{f(\bar{z})}}$.\ Obviously, that\ $f^{\ast}\in\bep$\ too.    An
entire function\ $f$\ is said to be   {\it real function} if\  $f\equiv f^{\ast}$.\  A real function\ $f\in\bep$\  takes
 only real values on\ $\R$\ , and every its complex zero necessarily is  conjugate zero.

\ \ \ \ \ If\ $f$\ is  an entire function of exponential type, then
\[
h_{f}(\theta)=\limsup_{r\to\infty}\frac{\ln |f(r e^{i\theta
})|}{r},\quad \theta\in [0;2\pi],
\]
is called the {\it indicator function} for\ $f$\ (with respect the first order of growth), and
\[
 \sigma_f=\max_{\theta\in[0;2\pi]}|h_{f}(\theta)|.
 \]
 An entire function\ $F$\ of exponential type is called  $\sigma$-{\it sine-type function} (or simple
sine-type function),  if there are  positive constants\ $c_1, c_2$,\ and\ $K$\ such that
\begin{equation}
c_1\le |F(x+iy)| e^{-\sigma|y|}\le c_2,\quad x, y\in\R, \  |y|\ge K.
\end{equation}
The notion of sine-type functions was introduced by  B.Ya. Levin in\ [5],\ where another but equivalent definition
was indicated. These functions compose the wide class,  and by various  authors the numerous applications of this class
were  found. For example, it contain any function of the form
\[
F(z)=\int_{-\sigma}^{\sigma}e^{-itz}d\mu(t),
\]
where\ $\mu$\ is any finite complex  measure such that\ $\mu(\{-\sigma\})\neq 0$,\ and $\mu(\{\sigma\})\neq 0$.\
Finally, every $\sigma$-sine-type function\ $F$\  has the type\ $\sigma_F=\sigma$\ and belongs to the Bernstein
space\ $B^{\infty}_{\sigma}$.

\ \ \ \ \ The exposed points of\ $\dbe$\ have the following characterization in terms of entire functions (see\ [2]).

{\bf Theorem 2.1.}\quad  \textit{
Let\ $f\in\be$\ and\ $\|f\|=1$.\ Then}\ $f\in\expp$\ \textit{if and only if:  \newline
(i)\ $\sigma_f=\sigma$; \quad (ii)\ $f$\ has no conjugate zeros; \quad
(iii)\ the real zeros of\ $f$\  are all simple; \qquad
(iv)\ $\int_{-\infty}^{\infty}|f(x)|h(x)
dx=\infty$\ whenever\ $h$\ is an entire function such that\ $\sigma_h=0$\ and}\ $h\ge 0$\ on\ $\R$.

\ \ \ \ \ Now we present the predual Banach space to \ $\be$.\ To this end, we consider the duality pair\
$(C_0(\R), M(\R))$,\ where\
$C_0(\R)$ is the usual Banach space of complex continuous functions on\ $\R$\ vanishing at infinity, and\ $M(\R)$\ is
the Banach algebra of all finite regular complex Borel measures on\ $\R$\ equipped with the total  variation norm,\
and convolution as multiplication. Let\ ${\mathfrak I}_\sigma$\ be  the
closed ideal of  all\  $\mu\in M(\R)$\ such that the Fourier--Stieltjes transforms\
${\hat \mu}(t)$\ vanish  for\ $|t|\ge \sigma$.\ Set
\begin{equation}
C_{0,\sigma}(\R)=\big\{ f\in C_0(\R)\colon\   \int\limits_{\R} f(x) \dr \mu(x)=0,\ \forall\mu\in
{{\mathfrak I}}_\sigma\big\}.\
\end{equation}
{\bf Proposition 2.2.} \ [6]\quad  \textit{ The space\ $\be$\ is the dual space of \ $C_0(\R)/ C_{0,\sigma}(\R)$.}

\vspace{4mm}
\section{ EXPOSED AND STRONGLY EXPOSED FUNCTIONS  }

\ \ \ \ \ \ Let\ $\Phi$\ be a continuous linear functional on\ $\be$,\ and suppose\ $\Phi$\ attains its  norm.
We  call a function\  $f\in\dbe$\ with\ $\|f\|=1$\ an {\it extremal} of\ $\Phi$\ if\ $\Phi(f)=\|\Phi\|$.\
   It may be noted that\  $f\in{\text{ exp}}\dbe$\  if and only if there exists\
$\Phi\in(B^1_{\sigma})^{\ast}$\ such that\ $\Phi$\  has an unique  extremal\ $f$.\  On the other hand, by the
Hahn-Banach theorem, every\ $f\in\be$,\  $\|f\|=1$,\ is an extremal for some\
$\Phi\in (B^1_{\sigma})^{\ast}$.\  We select among such functionals the following
\begin{equation}
\Phi_f(g)=\int_{\R} g(x) u_f(x) dx,\qquad g\in\be,
\end{equation}
where\ $u_f (x)$\ is  the function\ $f^{\ast}(x)/|f(x)|\in L^{\infty}(\mathbb{R})$\ defined for almost
all\ $x\in\R$.\ Thus, if\ $f\in{\text{ exp}}\dbe$,\ then\ $\Phi_f$\ is an exposing functional for\ $f$.\

\ \ \  \ \ Note that every\ $f\in{\text{ exp}}\dbe$\  has a unique exposing functional, which coincides
with\ (3.1).\   Indeed, assume that\
$\Phi\in(\be)^{\ast}$\
expose\ $f\in{\text{ exp}}\dbe$.\  By the Hahn--Banach  theorem,\ $\Phi$\   can be continued up to a linear functional\
$\Psi$\  on\ $L^1(\mathbb{R})$\  without increase of its norm. Then there is\ $\psi\in L^{\infty}(\mathbb{R})$,\
 $\|\psi\|=1$,\ such that\ $\Psi(a)=\int_{\R}a(t){\overline \psi(t)}dt$\ for all\ $a\in L^1(\mathbb{R})$.\
Now since\ $\Psi(f)=\Phi(f)= 1$,\  we obtain that\ $\psi(t)$\ coincides with\ $f(t)/|f(t)|$\ for almost
all\ $t\in\R$.\ Therefore,\ $\Phi=\Phi_f$.

\ \ \ \ \ The following  theorem gives simple and easily verified sufficient conditions of exposedness in\ $\dbe$.\

\begin{teo}
Let\ $f\in\be$,\  $\|f\|=1$,\ and\ $\sigma_f=\sigma$.\ Suppose\ $f$\  has no conjugate zeros and all real zeros of\ $f$\
are simple. If there exist\ $\tau$,\ $0<\tau\le 3$\ and\ $y_0\in\R$\  such that
\begin{equation}
\inf_{x\in\R}(|x+iy_0|^{\tau}|f(x+iy_0)|)>0,
\end{equation}
 then\  $f\in\expp$.\
\end{teo}

{\bf Remark 3.2.}\quad  The Plancherel--Polya inequality\ (2.1)\ implies that each\ $f\in\be$\  belongs to\
$L^1(\R)$\ not only on\ $\R$,\ but also on each line\ $\R+ia=\{z\in\Co: z=x+ia, x\in\R\}$,\  $a\in\R$.\
Therefore,\ $f_a(z):=f(z+ia)$\ belongs to\ $\be$\  for any\ $a\in\R$.\  Thus, by\ (2.2),\ we conclude that
\begin{equation}
\lim_{x\to\pm\infty} f(x+ia)=0, \quad a\in\R.
\end{equation}
\ \ \ \ \ \ {\bf Proof of Theorem 3.1.}\quad Assume\ $h$\ is a nonnegative entire function such that\ $\sigma_h=0$\
and\ $hf\in\be$.\ In view of Theorem 2.1, it suffices to show that\ $h$\ is a constant.  Set\
$g:=h f$.\ Since\ $g\in\be$,\ then combining\ (3.2)\ and\ (3.3),\ we see that there exist\ $\tau\in (0, 3]$\ and\
$y_0\in\R$\  such that
\begin{equation}
|h(x+iy_0)|= o (|x|^{\tau}),\quad x\to\pm\infty.
\end{equation}
The Phragmen-Lindel\"{o}f theorem implies that  entire function\ $h$\  of minimal type\ $\sigma_{h}=0$\ must be
a polynomial. Moreover, since\ $\tau\le 3$,\ then\ (3.4)\ implies  that the degree of this polynomial
satisfies\ $\deg h\le 2$.\  Because\ $h(x)\ge 0$,\  $x\in\R$,\ it follows that either\
$\deg h=0$\  or\ $\deg h=2$.\ It remains only to show that\ $\deg h=0$.\
To this end, let us suppose on the contrary, that\ $\deg h=2$.\ Since\ $h f\in\be$,\ then\
$\varphi(z) :=z^2 f(z)$\ also belongs to\ $\be$.\  On the other hand, if
\[
a= \inf_{x\in\R}(|x+iy_0|^{\tau}|f(x+iy_0)|),
\]
then\ (3.2)\ implies that
\[
|\varphi(x+iy_0)|\ge a |x+iy_0|^{2-\tau},\quad x\in\R.
\]
 By\ (3.2),\ we have\ $a>0$.\ As\ $\tau\le 3$,\ we see that the trace of\ $\varphi$\
 on the line\ $\R+iy_0$\ does not belongs to the space\ $L^1$.\ Hence\
$\varphi\not\in\be$,\ by Remark 3.4. This contradiction completes the proof.

\ \ \ \ \ \ Consider relations between sine-type functions  and exposedness.  We show  that each
such function determines a large set of exposed functions in\ $\dbe$.\ Before stating the next proposition,
we note that each sine-type function\ $F$,\  $F(z)\not\equiv ce^{iaz}$,\ $c\in\Co$,\ $a\in \R$,\  has infinitely
many zeros\ (see\ [4]\ for  details). Then for any\ $n\in\mathbb{N}$\ there exists a polynomial\ $p$\ such that\
$\deg p=n$\ and\ $N(p)\subset N(F)$.\

{\bf Proposition 3.3.}\quad \textit{
Let\ $F$\ be a\  $\sigma$-sine-type function,\ $F(z)\not\equiv ce^{\pm i\sigma z}$,\ $c\in\Co$.\ Let\ $F$\ has neither
conjugate nor  multiple real zeros. Suppose\ $q$\ is a polynomial such that\ $N(q)\subset N(F)$.\ Set
\begin{equation}
 f_q(z)=\alpha \frac{F(z)}{q(z)},
 \end{equation}
where\ $\alpha$\ is any complex number such that\ $\|f_q\|_{L^1}=1$.\
Then\ $f_q$\ is an exposed point of\ $\dbe$\ if and only if \ $2\le\deg q\le 3$.}

\ \ \ \ \ \ {\bf Proof.}\quad  Suppose\ $2\le \deg q\le 3$.\ There exists an\ $M>0$\ such that\
$N(q)\subset \{z\in\Co: |Im\, z|\le M\}$.\ Hence it follows from\ (2.3)\ and\ (3.5)\ that  if\ $y_0\in\R$\
and\ $|y_0|\ge \max (K; M)$,\ then
\[
|x+iy_0|^{\deg q} |f_q(x+iy_0)|\ge \alpha c_1e^{\sigma|y_0|}\frac{|x+iy_0|^{\deg q}}{|q(x+iy_0)|},
\]
where\ $c_1>0$.\ An application of Theorem 3.1 shows that\ $f_q\in\expp$.

\ \ \ \ \ If\ $\deg q <2$,\ then\ (2.3)\ implies that\ $f_q\notin L^1(\R)$.\ Thus\ $f_q\notin \be$.\ Finally, assume\
$\deg q> 3$.\ Then the function\ $z^2f_q(z)$\ belongs to\ $\be$,\ i.e.,\ $f_q$\ has entire and
nonnegative on\ $\R$\ multiplier\ $h(z)=z^2$\ with\ $\sigma_h=0$.\ By Theorem 2.1, we get  that\ $f_q\notin\expp$.

\ \ \ \ \ \ {\bf Example 3.4.}\quad There exists an exposed function\ $f$\  in\ $\dbe$\ such that its zeros set \ $N(f)$\
is not separated.  Recall that\ $\Lambda\subset \mathbb{C}$\  is called {\it separated}, if

\[
\inf_{\substack{z, \lambda\in \Lambda \\z\neq \lambda}} |z-\lambda|= \delta>0.
\]
Let\ $0<\varepsilon<\pi/2$\ and
\begin{equation}
 f(z)=\alpha\frac{\cos( \sigma\frac{z}2) \cos\sqrt{(\sigma\frac{z}2)^2+\varepsilon^2}}{\sigma^2z^2-\pi^2},
 \end{equation}
where\  $\alpha$\ is  is any complex number such that\ $\|f\|_{L^1(\R)}=1$.\
 Since\ $0<\varepsilon<\pi/2$,\ then the $\sigma$-sine-type function
 \[
F(z)=\cos\Bigl( \sigma\frac{z}2\Bigr) \cos\sqrt{\Bigl(\sigma\frac{z}2\Bigr)^2+\varepsilon^2}
\]
has only real and simple zeros. Therefore\ $f\in \expp$,\ by Proposition 3.3, and
\[
N(f)=\Bigl\{\frac2{\sigma}\bigl(\frac{\pi}2+\pi k\Bigr), \ k\in{\mathbb{Z}}\Bigr\}\cup\Bigl\{\pm\frac2{\sigma}
\sqrt{\Bigl(\frac{\pi}2+\pi l\Bigr)^2-\varepsilon^2}, \ l\in{\mathbb{Z}}\Bigr\}
\]
 is obviously not separated.

{\bf Proposition 3.5.}\quad \textit{Let\ $f\in\dbe$,\ and assume\ $N(f)\cap\R\neq\emptyset$.\ If}\ $f\in \exps$,\
  \textit{  then\ $N(f)\cap\R$\   is separated set.}

\ \ \ \ \ In view of this proposition, if\ $0<\varepsilon<\pi/2$,\  then the function\ (3.6)\ is an example
of exposed point of\ $\dbe$\ such that it is not strongly exposed. Thus we have:

{\bf Corollary  3.6.}
\[
\exps\subsetneqq\expp.
\]
\ \ \ \ \ \ {\bf Proof of Proposition 3.5.}\quad Let\ $f\in\expp$.\ Suppose\ $f$\ has not separated set of real zeros,
i.e. there exists a double sequence\ $(x_n, y_n)$,\ $n\in \mathbb{N}$,\ where\ $x_n, y_n\in N(f)\cap\R$,\ such that\
$0<y_n-x_n<1/n$.\ Let's choose arbitrarily\
$a\in\R$ such that
 \begin{equation}
f(a)\not= 0.
 \end{equation}
Put
\[
f_n(x)=\alpha_n(x-a)^2\frac{f(x)}{(x-x_n)(x-y_n)},\quad n\in\mathbb{N},
\]
where\ $\alpha_n$\ is such positive numbers that\ $\|f_n\|=1$.\ If\ $\Phi_f$\ is the exposing functional\ (3.1)\
for\ $f$,\ then
 \begin{equation}
\Phi_f(f_n)=\int_{-\infty}^{\infty} |f_n(x)| dx-2\int_{x_n}^{y_n} |f_n(x)| dx= 1- 2\int_{x_n}^{y_n} |f_n(x)| dx.
\end{equation}
It is known that if\ $g\in\be$\ and\ $u, v\in\R$,\ $u\le v$,\ then
\[
\int_u^v |g(x)| dx\le A(v-u)\|g\|,
\]
where\ $A>0$\ do not depends from $u, v$\ and\ $g$\ (the exact value of\ $A$\ is obtained in\ [7]).\  Combining this
with\ (3.8),\ we obtain
  \begin{equation}
\lim_{n\to\infty} \Phi_f(f_n)=1.\
 \end{equation}
\ \ \ \ \ By definition of\ $f_n$,\ we have\ $\lim_{n\to\infty} f_n(a)=0$.\ Then\ (3.7)\ implies, in particular,
that\ $(f_n)$\  does not converges  to\ $f$\ uniform on compact
sets in\ $\R$.\ Hence (see\ [3]) we conclude\ $f_n\nrightarrow f$\ in\ $\be$.\   Combining this with\ (3.9),\ we see
that\ $f$\ is not strongly exposed\ in\ $\dbe$.
\vspace{4mm}
\section{ WEAK$^{\ast}$ EXPOSED FUNCTIONS  }

\ \ \ \ \ Since\ $\be$\ is  a closed subspace of\ $L^1(\R)$,\  the dual space of\ $\be$\ is  isometrically isomorphic
to  \ $L^{\infty}(\R)/{\text{Ann}}( \be)$,\ where\
${\text{Ann}}( \be)$\  is the annihilator of\ $\be$, i.e.
\[
{\text{Ann}}( \be)=\Bigl\{g\in L^{\infty}(\R): \ \int_{-\infty}^{\infty}g(t){\overline{f(t)}} dt=0,\ f\in\be\Bigr\}.
\]
 ${\text{Ann}}( \be)$\  has also a direct  description in terms of  spectrum\ ${\text{spec}} f$\  of a
function\ $f\in L^{\infty}(\R)$.\  By\ ${\text{spec}} f$\ we mean the
support of Fourier transform of\ $f$ (understood in the distributional sense). On the other hand, it is easy to see
that\ ${\text{Ann}}( \be)$\  coincides with a
subspace in\ $L^{\infty}(\R)$\ of {\it high pas} functions or functions  in\
$L^{\infty}(\R)$\ with the spectral gap\ $(-\sigma,\sigma)$:
\[
{\text{Ann}}( \be)= H^{\infty}_{\sigma}(\R) :=\{ f\in L^{\infty}(\R):\ {\text{spec}}
f\cap (-\sigma,\sigma)=\emptyset\},
\]
(see also\ [9]).\ Hence each\ $\Psi\in(\be)^{\ast}$\ can be set as
 \begin{equation}
\Psi(g)=\int^{\infty}_{-\infty} g(t){\overline{\psi(t)}} dt,
 \end{equation}
 where\ $\psi\in L^{\infty}(\R)$.\  In particular, this is not a unique such representation of\ $\Psi$.\ We can
 replace\  $\psi$\ in\ (4.1)\  by\  $\psi+h$,\ where\  $h$\ is any function in\  $\he$.\

 \ \ \ \ \ \ We present now some properties of weak$^{\ast}$\ continuous functionals on\ $\be$.\  The definition\ (2.4)\
 of\ $C_{0,\sigma}(\R)$\ implies that\ $C_{0,\sigma}(\R)\subset \he$.\  Combining  this with expression\ (4.1),\ we
 obtain:

 {\bf Proposition 4.1.}\quad \textit{A  functional defined by}\ (4.1)\ \textit{is weak$^{\ast}$\ continuous on\ $\be$\
 if and  only if there  are\ $u_{\psi}\in C_{0}(\R)$\  and\ $h_{\psi}\in\he$\ such that\ $\psi=u_{\psi}+h_{\psi}$.}

 {\bf Lemma 4.2.}\quad \textit{If\ $\psi\in L^1(\R)\cap L^{\infty}(\R)$,\ then}\ (4.1)\
  \textit{is weak$^{\ast}$ continuous functional on\ $\be$ .}

\ \ \ \ \ \ {\bf Proof.}\quad  The function\ $S_{\sigma}(x):=\sin \sigma x/\pi x$\ belongs to\
$B^p_{\sigma}$,\  for every\ $1< p\le\infty$.\ Moreover, the Fourier transform of\ $S_{\sigma}$\ coincides with the
 characteristic (indicator)  function of\ $[-\sigma,\sigma]$.\ Hence the  operator\ $f\to S_{\sigma}$\ on\ $\be$\
is well defined  and  coincides with the identity operator. Therefore, we have
\begin{equation}
f(t)=\int_{-\infty}^{\infty} f(x) S_{\sigma}(t-x) dx,\quad f\in\be.
\end{equation}

\ \ \ \ \ If\ $\psi\in L^1(\R)$,\ then the integral
\begin{equation}
\psi_1(x)=\int_{-\infty}^{\infty} \psi(x) S_{\sigma}(t-x) dt
\end{equation}
converges absolutely, and\ $\psi_1\in L^p(\R)$\ for every\ $1<p\le \infty$,\  by Young inequality.
    Assume now that\ $\psi\in L^1(\R)\cap L^{\infty}(\R)$.\ Then\ $\psi_1\in L^r(\R)$,\ for every\
$1\le p\le \infty$.\ Using the Parseval formula in\ $L^2(\R)$\ and applying\ to\ (4.3)\ the  Paley-Wiener theorem, we
see that\ $\psi_1\in B^2_{\sigma}$.\ By\ (2.2),\ this means that\ $\psi_1\in C_0(\R)$.\ Finally, we
will show  that functions\ $\psi$\ and\ $\psi_1$\ define in\ (4.1)\ the same functional on\ $\be$.\ This is equivalent
to stating that  $\psi-\psi_1\in H^{\infty}_{\sigma}(\R)$.\ Since\ (4.3)\ converges absolutely, then it follows
from\    (4.3)\ and Fubini's theorem that
\begin{eqnarray*}
\int^{\infty}_{-\infty}f(x){\overline{\psi_1(x)}}dx =\int^{\infty}_{-\infty}f(x)\Bigl[\int^{\infty}_{-\infty}
{\overline{\psi(t)}}S_{\sigma}(x-t) dt \Bigr]dx \\
 =\int^{\infty}_{-\infty}\Bigl[\int^{\infty}_{-\infty}f(x)  S_{\sigma}(t-x) dx\Bigr]
 {\overline{\psi(t)}}dt=\int^{\infty}_{-\infty}f(t){\overline{\psi(t)}}dt.
\end{eqnarray*}
This concludes the proof.

{\bf Remark 4.3.}\quad The statement of Lemma 4.2 carries over to\ $\psi\in L^p(\R)\cap L^{\infty}(\R)$,\ for
each\ $1\le p<\infty$.

 {\bf Theorem 4.4.}\quad \textit{Any exposed$^{\ast}$\ function in\ $\dbe$\ is strongly exposed$^{\ast}$,\ i.e.}
\begin{equation}
\exppp=\expss.
\end{equation}
\ \ \ \ \ \ {\bf Proof.}\quad Let\ $f\in \exppp$,\ and let\ $\Phi_f$\ be an exposing functional for\ $f$.\ Suppose\
$(f_n)$\ is a sequence in\ $\dbe$\ such that\ $\lim_{n} \Phi_f(f_n)=1$.\ We may, without loss of generality, assume
that\ $\|f_n\|=1$,\  $n\in\mathbb{N}$.\ We will show that\ $(f_n)$\ converges to\ $f$\ in\ $\be$.\ This  is equivalent to
the following statement: given any subsequence\ $(f_{n_{k}})_k$\ of\ $(f_n)$\  there exists  convergent to\ $f$\
subsequence\ $(f_{n_{k_{j}}})_j$\ of\ $(f_{n_{k}})_k$.\ Now fix a subsequence\ $(f_{n_{k}})_k$.\ Put\ $F_k =f_{n_{k}}$,\
$k\in \mathbb{N}$,
  to simplify the writing.  \ $\dbe$\ is compact set with respect to the topology of uniform convergence on
  compact subsets of\ $\R$\ (see\ [3]).\        Thus\ $(F_k) $\ has a subsequence\ $(F_{k_{j}})_j$\ such that\
  $(F_{k_{j}})_j$\
  converges to  some\ $g\in\dbe$\ uniformly  on   compact subsets of\ $\R$.\ Since\ $f\in \exppp$\ and each exposed
  function of\ $\dbe$\ has an unique  exposing functional\ $\Phi_f$,\ then\ $\Phi_f$\ is weak$^{\ast}$\ continuous.
  Therefore\ 1=$\lim_{j}\Phi_f(F_{k_{j}})=\Phi_f(g)$,\  which gives that\ $g=f$,\ because\ $f$\ is exposed. Finally, we
  will show that\ $(F_{k_{j}})$\ converges to\ $f$\ in\ $\be$.\ Let\ $\varepsilon>0$.\ Then there exists a compact
  subset\ $\Omega\subset\R$\ such that
\begin{equation}
\int_\Omega |f(t)|dt> 1-\varepsilon.
\end{equation}
Now fix\ $\Omega$\ and choose\ $M\in \mathbb{N}$\ such that
\begin{equation}
\max_{t\in\Omega} |F_{k_{j}}(t) - f(t)| < \varepsilon/2|\Omega|,\quad j> M,
\end{equation}
where\ $|\Omega|$\ is Lebesgue measure of\ $\Omega$.\ Since\ $\|f\|=1$,\ it follows from\ (4.5)\ that \
 $\int_{\Omega}|F_{k_{j}}|\ge 1-2\varepsilon$. \ Combining this with\ (4.5)\ and\ (4.6),\ and using\
 $\|F_{k_{j}}\|=1$,\ $j\in \mathbb{N}$,\ we obtain
\[
\|F_{k_{j}}- f\|\le \int_{\Omega}|F_{k_{j}}(t)-f(t)|dt +\int_{\R\setminus\Omega}|F_{k_{j}}(t)|dt+
\int_{\R\setminus\Omega}|f(t)|dt <4\varepsilon,\quad j> M.
\]
Therefore,\ $(F_{k_{j}})$\ converges to\ $f$\ in\ $\be$,\ and\ $f\in \expss$.\ This concludes the proof.

\ \ \ \ \ We need the following  technically  lemma.

{\bf Lemma 4.5.}\quad \textit{Let\ $f$\ be an exposed function in\ $\dbe$\ with exposed functional\ $\Phi_{f}$.\
 Assume\ $f_n \in \dbe$,\ $n\in \mathbb{N}$,\ and\ $\lim_{n}\Phi_{f}(f_n)=1$.\ If\ $(f_n)$\ converges uniformly on compact
sets in\ $\R$\ to\ $f_0\in\dbe$,\ then there exists\ $\alpha$,\ $0\le \alpha\le 1$,\ such that\ $f_0=\alpha f$.}

\ \ \ \ \ \ {\bf Proof.}\quad Let\ $\varepsilon>0$.\ There exists a compact subset\ $K\subset\R$\ such that
\begin{equation}
\int_{K}|f_0(x)| dx\ge \|f_0\|-\varepsilon.
\end{equation}
Let\ $u_f$\ is the unimodular function\ $f^{\ast}/|f|$,\ defined on\ $\R\setminus N(f)$.\ Then\ (3.1)\  and\ (4.7)\
imply that
\begin{equation}
\Bigl|\int_{\R} f_0(x) u_f(x) dx\Bigr|\ge \Bigl|\int_{K} f_0(x) u_f(x) dx\Bigr|-\varepsilon.
\end{equation}
Choose\ $M_1>0$\  so that
\[
\max_{t\in K}| f_n(x)- f_0(x)|\le \frac{\varepsilon}{|K|},\quad n>M_1.
\]
 Thus if\ $n>M_1$,\ then
\begin{equation}
\int_{K}|f_n(x)| dx > \int_{K}|f_0(x)| dx -\varepsilon,
\end{equation}
and
\begin{equation}
\Bigl|\int_{\R} f_0(x) u_f(x) dx\Bigr|> \Bigl|\int_{K} f_n(x) u_f(x) dx\Bigr|-\varepsilon.
\end{equation}
Similarly  there exists an\ $M_2>0$\  such that
\begin{equation}
\Bigl|\int_{K} f_n(x) u_f(x) dx\Bigr|> \int_{K} |f_n(x)| dx-\varepsilon, \quad  n>M_2.
\end{equation}
Actually, assume, to the contrary, that there exists a subsequence\ $(f_{n_{k}})$\ such that
\[
\Bigl|\int_{K} f_{n_{k}}(x) u_f(x) dx\Bigr| \le \int_{K} |f_{n_{k}}(x)| dx-\varepsilon,
\quad k\in\mathbb{N}.
\]
Then we conclude that
\[
|\Phi_f(f_{n_{k}})|\le \Bigl|\int_{K} f_{n_{k}}(x) u_f(x) dx\Bigr| +
 \Bigl|\int_{\R\setminus K}  f_{n_{k}}(x) u_f(x) dx\Bigr|\le
 \]
 \[
 \int_{K} |f_{n_{k}}(x)| dx-\varepsilon+
 \int_{\R\setminus K} |f_n(x)| dx\le 1-\varepsilon,
 \]
contrary to\ $\lim_k \Phi_f(f_{n_{k}})= \lim_n \Phi_f(f_{n})=1$.\ Combining\ (4.7)\ --\ (4.11),\ we see that
\[
\Phi_f(f_0)|> \|f_0\|-5\varepsilon.
\]
Therefore, if\ $f_0\not\equiv 0$,\ then there exists \ $c\in \mathbb{C}$,\ $|c|\ge 1$,\ such that\ $\|c f_0\|=1$\ and\
$c f_0$\ is an extremal of\ $\Phi_f$.\ Since $\Phi_f$\ exposes\ $f\in\expp$,\ then\ $f_0=\alpha f$,\ where\
$\alpha=1/c$ .

\ \ \ \ \ To complete the proof, we show that\ $\alpha$\ is a positive number. As was proved above,  there exists
a sequence\ $(\varepsilon_n)$,\ $\varepsilon_n\ge 0$,\ $n\in \mathbb{N}$,\ such that
\begin{equation}
\|f_n -f_0\|=\|f_n-\alpha f\|\le 1-|\alpha|+\varepsilon_n,\quad n\in \mathbb{N}.
\end{equation}
Suppose\ $\alpha\neq 0$,\ and put
\[
h_n(x):= \frac{f_n(x)-f_0(x)}{1-\alpha},\quad n\in \mathbb{N}.
\]
Then\ $\lim_{n}\Phi_f(h_n)=1$,\ and  hence\ $\|h_n\|\le 1$,\ $n\in \mathbb{N}$.\ Combining this information with\  (4.12),\
we see that\ $\alpha\in\R$\ and\ $0\le \alpha\le 1$.

{\bf Theorem 4.6.}\quad \textit{Let
\begin{equation}
f(x)=\alpha\frac{2\cos\sigma x-1}{(3\sigma x)^2-\pi^2},
\end{equation}
where\ $\alpha>0$ is such  that\ $\|f_a\|=1$.\ The function}\ (4.13)\ \textit{is  strongly exposed in\ $\dbe$, but
 not weak$^{\ast}$ exposed.\ In particular, }
\[
\exppp\subsetneqq \exps.
\]
\ \ \ \ \ \ {\bf Proof.}\quad Suppose\ $F(x)=2\cos\sigma x-1$,\ and\ $q(x)=(3\sigma x)^2-\pi^2$.\ Since\ $F$\ is a
$\sigma$-sine-type function, it follows from
Proposition 3.3 that the function\ (4.13)\ belongs to\ $\expp$.\

\ \ \ \ \ We now show that\ $f\not\in\exppp$.\ If\ $x\in\Omega:=
\R\setminus\bigl\{\pm \pi/3\sigma +2\pi k/\sigma, \ k\in\mathbb{Z}\bigr\}$,\ then
\begin{equation}
\frac{f^{\ast}(x)}{|f(x)|} = -2\chi_{\sigma}(x)+\frac{F^{\ast}(x)}{|F(x)|},
\end{equation}
where\ $\chi_{\sigma}$\ is the characteristic function of\ $[-\pi/3\sigma, \pi/3\sigma]$.\ The unimodular function\
$F^{\ast}/|F|$\ is periodic with period\ $2\pi/\sigma$,\ and has finite variation in\ $[0, 2\pi/\sigma]$.\ Hence\
$F^{\ast}/|F|$\ coincides on\ $\Omega$\ with its Fourier series
\begin{equation}
\frac{F^{\ast}(x)}{|F(x)|} = -\frac13 +\frac{4}{\pi}\sum_{n=1}^{\infty} \frac{\sin\bigl(\frac{\pi n}3\bigr)}{n}
\cos( \sigma n x),\quad x\in \Omega.
\end{equation}
Put
\[
h_{\sigma}(x)=\frac{4}{\pi}\sum_{n=1}^{\infty} \frac{\sin\bigl(\frac{\pi n}3\bigr)}{n}
\cos( \sigma n x).
\]
Then\ $h_{\sigma}\in L^{\infty}(\R)$.\ The Fourier transform of\ $h_{\sigma}$\ (in the distributional sense)
coincides with
\[
\sum_{n=1}^{\infty}\frac{\sin\bigl(\frac{\pi n}3\bigr)}{n}\bigl(\delta(-\sigma n)+\delta(\sigma n)\bigr),
\]
where\ $\delta(t)$\ is the Dirac distribution with support at\ $t\in\R$.\ Hence\ $h_{\sigma}\in H^{\infty}_{\sigma}(\R)$.\
Combining\ (4.14)\ and\ (4.15),\ we find that exposing functional\ (3.1)\ of\ $f$\ admits  the following representation
\begin{equation}
\Phi_f(g)=I_f(g)+K_f(g),\quad g\in\be,
\end{equation}
where
\[
K_f(g)=\int_{-\infty}^{\infty}\bigl(h_{\sigma}(x)-2 \chi_{\sigma}(x)\bigr) g(x) dx,
\]
and
\[
I_f(g)=-\frac13\int_{-\infty}^{\infty}g(x) dx.
\]
It follows from Proposition 4.1  and Lemma 4.2  that\ $K_f$\ is weak$^{\ast}$ continuous on\ $\be$.\ On the other hand,
\ the functional $I_f$\ is not weak$^{\ast}$ continuous. Actually every weak$^{\ast}$ continuous functional on\ $\be$\
is also continuous under the topology of\ $\be$-norm bounded  uniform convergence on compact subsets of\ $\R$.\ If
\[
g(x)=\biggl(\frac{\sin\frac{\sigma x}{2}}{x}\biggr)^2,
\]
then\ $g_n(x)= g(x-n)$, $n\in \mathbb{N}$,\ is  bounded in\ $\be$-norm and converges uniformly on compact subsets
of\ $\R$\ to zero function. However, the following is true:
\[
\lim_{n\to\infty} I_f(g_n)= -\frac13\|g_n\|=-\frac13\|g\|<0.
\]
It follows that\ $I_f$,\ and hence that\ $\Phi_f$\ are not weak$^{\ast}$ continuous  on\ $\be$.\ Thus, we obtain that\
$f\not\in\exppp$.\

\ \ \ \ \ Finally, we show that\ $f\in\exps$.\ This is the same as: given any subsequence\ $(f_{n_{k}})_k$\ of\ $(f_n)$\
there exists  convergent to\ $f$\ subsequence\ $(f_{n_{k_{j}}})_j$\ of\ $(f_{n_{k}})_k$.\ Put\ $F_k =f_{n_{k}}$,\
$k\in \mathbb{N}$,\   to simplify the writing. Thus\ $(F_k) $\ has an subsequence\ $(F_{k_{j}})_j$\ such that\
$(F_{k_{j}})_j$\  converges to  some\ $f_0\in\dbe$\ uniformly  on   compact subsets of\ $\R$.\ By Lemma 4.5, there
is\ $\alpha$,\   $0\le\alpha\le 1$,\ such that\ $f_0=\alpha  f$.\ We show that\ $\alpha=1$.\ To this end, assume, on
the contrary,   that\ $0\le \alpha<1$.\  Put
\[
g_j=\frac1{1-\alpha}\Bigl(F_{k_{j}}-\alpha f\Bigr),\quad j\in\mathbb{N}.
\]
Then\ $g_j\in\dbe$,\ $\lim_{j} \Phi_f(g_j)=1$,\ and\ $(g_j)$\ converges uniformly on compact subsets
of\ $\R$\ to zero function. Thus\ $\lim_{j} K_f(g_j)=0$,\ and\ (4.13)\ implies
\[
\lim_{j\to\infty}|\Phi_f(g_j)|\le \lim_{j\to\infty}|I_f(g_j)|=\frac13.
\]
This contradicts the assumption that\ $0\le \alpha<1$.\ Hence\ $f_0=f$,\ and this proves that\ $f\in\exps$.

\vspace{10mm}

 \centerline{\bf{REFERENCES}}
 \vspace{4mm}

1. \ E. Asplund,\ Fr\'{e}chet differentiability of convex functions.  Acta Math. {\bf 121}, 31-47 (1968).

2. \  K.M. Dyakonov,\ Polynomials and entire functions: Zeros and geometry of the unit ball. Math. Res. Lett. {\bf 7},
 No.4, 393-404 (2000).

3. \ J.R. Higgins,\ Sampling theory in Fourier and Signal analysis: Foundations. Clarendon Press, Oxford 1996.

4. \  B.Ya.  Levin,\ {\it Lectures on Entire Functions},\ Translations of Mathematical Monographs, {\bf 150}, American
Math. Soc. (1996).

5. \ B.Ya. Levin,\ On bases of exponential functions in\ $L^2$.\ Zap. Mekh. Mat. Fak. Kharkov Gos. Univ. i
Kharkov Mat. Obshch. {\bf 27}(4), 39-48 (1961). (Russian).

6. \ S. Norvidas,\ Majorants and extreme points of unit balls in Bernstein spaces. Lith. Math. J. {\bf 44}(1),
78-84  (2004).

7. \ S. Norvidas,\ On localization of functions in the Bernstein space. Lith. Math. J. {\bf 47}(4), 470-483  (2007).

8. \ R.R. Phelps,\ Dentability and extreme points in Banach space. J. Funct. Anal. {\bf 16}, 78-90 (1974).

9. \ H.S. Shapiro,\ Topics in approximation theory. Lecture Notes in Mathematics. 187. Berlin-Heidelberg-New York:
Springer-Verlag  1971.

}}
\end{document}